\documentclass[12pt,a4paper,reqno]{amsart}

\usepackage{amssymb}
\usepackage{amscd}
\usepackage{amsfonts}
\usepackage{setspace}
\usepackage{version}

\numberwithin{equation}{section}
\theoremstyle{plain}
\newtheorem{theorem}[subsection]{Theorem}
\newtheorem{proposition}[subsection]{Proposition}
\newtheorem{lemma}[subsection]{Lemma}
\newtheorem{corollary}[subsection]{Corollary}

\theoremstyle{definition}

\newtheorem{definition}[subsection]{Definition}

\newtheorem*{mainthm1-again}{Theorem \ref{mainthm1}'}
\newtheorem*{mainthm2-again}{Theorem \ref{mainthm2}'}
\newtheorem*{frei-ruz-again}{Theorem \ref{frei-ruz}'}
\renewcommand{\leq}{\leqslant}
\renewcommand{\geq}{\geqslant}
\newsavebox{\proofbox}
\savebox{\proofbox}{\begin{picture}(7,7)  \put(0,0){\framebox(7,7){}}\end{picture}}

\newcommand\C{\mathbb{C}}

\newcommand\GL{\operatorname{GL}}

\newcommand\id{\operatorname{id}}

\newcommand\Upp{\operatorname{Upp}}
\newcommand\Uni{\operatorname{Uni}}

\def\endproof{\hfill{\usebox{\proofbox}}\vspace{11pt}}

\onehalfspace
\begin{document}
\title{Approximate groups, II:  the solvable linear case}
\author{Emmanuel Breuillard}
\address{Laboratoire de Math\'{e}matiques, B\^{a}timent 425, Universit\'{e}
Paris-Sud, 91405 Orsay, France.}
\email{emmanuel.breuillard@math.u-psud.fr}
\author{Ben Green}
\address{Centre for Mathematical Sciences\\
Wilberforce Road\\
Cambridge CB3 0WA\\
England }
\email{b.j.green@dpmms.cam.ac.uk}
\subjclass{}

\begin{abstract}
We describe the structure of  ``$K$-approximate
subgroups'' of solvable subgroups of $\GL_{n}(\mathbb{C})$, showing that they have a large nilpotent piece. By combining this with
the main result of our recent paper on approximate subgroups of
torsion-free nilpotent groups \cite{breuillard-green-nil}, we show that such
approximate subgroups are efficiently controlled by nilpotent progressions.
\end{abstract}

\maketitle
\tableofcontents

\setcounter{tocdepth}{1}

\section{Introduction}

This paper is the second in a series concerning \emph{approximate groups}. The first paper in the series is \cite{breuillard-green-nil},
which discusses approximate subgroups of torsion-free nilpotent groups. The
reader is referred to that paper for a more extensive discussion on the
background to the material we discuss here. Let us recall the definition of
``$K$-approximate group''.

\begin{definition}[Approximate groups]
Let $G$ be some group and let $K \geq 1$. A set $A \subseteq G$ is called a $K$-approximate group if

\begin{enumerate}
\item It is symmetric, i.e. if $a \in A$ then $a^{-1} \in A$, and the
identity lies in $A$;

\item There is a symmetric subset $X$ lying in $A\cdot A$ with $|X| \leq K$
such that $A \cdot A \subseteq X \cdot A$.
\end{enumerate}
\end{definition}

Let us also recall Tao's notion of \emph{control} \cite{tao-solvable}.

\begin{definition}[Control]
Suppose that $A$ and $B$ are two sets in some ambient group, and that $K
\geq 1$ is a parameter. We say that $A$ is $K$-controlled by $B$, or that $B$
$K$-controls $A$, if $|B| \leq K|A|$ and there is some set $X$ in the
ambient group with $|X| \leq K$ and such that $A \subseteq (X \cdot B) \cap
(B \cdot X)$.
\end{definition}

\textsc{Solvable linear groups.} Our aim in this paper is to study $K$-approximate subgroups of solvable subgroups of $\GL_n(\mathbb{C})$,
for a fixed integer $n$. Recall that these are groups $G$ such that the
derived series 
\begin{equation*}
G_{(0)} = G, G_{(1)} = [G_{(0)},G_{(0)}], G_{(2)} = [G_{(1)},G_{(1)}],\dots
\end{equation*}
terminates with $G_{(s+1)} = \{\id\}$. For example the group $G =
\left( 
\begin{smallmatrix}
\mathbb{C}^{\times} & \mathbb{C} & \mathbb{C} \\ 
0 & \mathbb{C}^{\times} & \mathbb{C} \\ 
0 & 0 & \mathbb{C}^{\times}
\end{smallmatrix}
\right)$ is solvable, as is any subgroup of it.

This example is actually fairly typical, as the following well-known theorem
of Mal'cev \cite{malcev-solvable} (see also \cite[Theorem 3.6]{wehrfritz})
shows. In this result and the rest of the paper we write $\Upp_n(\mathbb{C})$ for the (solvable) subgroup of $\GL_n(\mathbb{C})$
consisting of upper-triangular matrices.

\begin{theorem}[Mal'cev]
Suppose that $G$ is a solvable subgroup of $\GL_n(\mathbb{C})$. Then
there is a normal subgroup $H \lhd G$ whose index in $G$ is bounded by some
function $F(n)$ of $n$ only and which is conjugate to a subgroup of $\Upp_n(\mathbb{C})$.
\end{theorem}

We cannot take $F(n)=1$, that is to say not all solvable subgroups are
conjugate to a group of upper triangular matrices. Indeed any non-abelian
finite group has an irreducible linear representation of dimension at least $2$.

We now state our main results.

\begin{theorem}
\label{mainthm2} Let $K\geq 1$ be a parameter, suppose that $G\subseteq 
\GL_{n}(\mathbb{C})$ is a solvable group, and suppose that $A\subseteq
G$ is a $K$-approximate group. Then $A$ is $K^{C_{n}}$-controlled by a $K^{C_{n}}$-approximate group $B\subseteq G$ which generates a nilpotent
group of step at most $n-1$ contained in a conjugate of $\Upp_{n}(\mathbb{C})$.
\end{theorem}

\emph{Remarks.} We recalled in \cite{breuillard-green-nil} that all
finitely-generated torsion-free nilpotent groups may be realised as linear
groups. The same is by no means true for solvable groups; indeed if a group $G$ embeds in $\Upp_n(\mathbb{C})$ then $[G,G]$ is nilpotent, a very
special property in the class of all solvable groups. The paper of Tao \cite{tao-solvable} applies to solvable groups in
general, and so is much broader in scope than our work. That paper comes
with extremely weak bounds, however, whereas we obtain polynomial dependence
on the approximation parameter $K$. The work of Sanders \cite{sanders} also applies to a more general class of group than ours. His bounds are in a sense better than Tao's, but his structural conclusions are somewhat weaker.

\begin{corollary}
\label{main-cor} Let $K \geq 1$ be a parameter, suppose that $G \subseteq 
\GL_n(\mathbb{C})$ is a solvable group, and suppose that $A \subseteq
G $ is a $K$-approximate group. Then $A$ is $e^{K^{C_n}}$-controlled by a
nilpotent progression of dimension at most $K^{C_n}$ and step at most $n-1$.
\end{corollary}

\emph{Remarks.} Let us make particular note at this point of our dependence
on Part I of the series \cite{breuillard-green-nil}, where
the concept of a nilpotent progression that we have invoked here is
introduced and developed at some length.

\section{Preliminaries from multiplicative combinatorics}

\label{basic-add-sec}

Let us begin by stating the basic properties of approximate groups. The next
proposition was also stated in \cite{breuillard-green-nil} as Proposition 2.1. It is
due to Tao \cite{tao-noncommutative}.

\begin{proposition}[Approximate groups and control]
\label{approx-gp-basics} Let $K \geq 1$ be a parameter and let $A$ be a set
in some ambient group $G$. If $n \geq 1$ is an integer we write $A^n = \{a_1
\dots a_n : a_1,\dots,a_n \in A\}$ and $A^{\pm n} = \{a_1^{\varepsilon_1}
\dots a_n^{\varepsilon_n} : a_1,\dots,a_n \in A,
\varepsilon_1,\dots,\varepsilon_n \in \{-1,1\}\}$.

\begin{enumerate}
\item If $\pi : G \rightarrow H$ is a homomorphism and if $A \subseteq G$ is
a $K$-approximate group then $\pi(A)$ is a $K$-approximate subgroup of $H$.

\item If $A$ is a $K$-approximate group then $|A^{\pm n}| = |A^n| \leq
K^{n-1}|A|$ and $A^n$ is $K^{n+1}$-controlled by $A$.

\item If $B,C$ are further subsets of $G$ and if $A$ is $K$-controlled by $B$
and $B$ is $K$-controlled by $C$, then $A$ is $K^2$-controlled by $C$.

\item If the doubling constant $|A^2|/|A|$ is at most $K$ then there is an $f_1(K)$-approximate group $B \subseteq A^{\pm 3}$ which $f_2(K)$-controls $A$
. If the tripling constant $|A^3|/|A|$ is at most $K$ then we may take $B =
A^{\pm 3}$.

\item If $A$ is a $K$-approximate group and if $A^{\prime }\subseteq A$ is a
subset with $|A^{\prime }| \geq |A|/K$ then $A^{\prime \pm 3}$ is an $f_3(K)$
-approximate group which $f_4(K)$-controls $A$. The same is in fact true
under the essentially weaker assumption that $|A^3| \leq K|A|$.
\end{enumerate}

All of the quantities $f_1(K),\dots,f_4(K)$ can be taken to be polynomial in 
$K$.
\end{proposition}

Let us recall also the following simple but surprisingly powerful
observation, which is essentially \cite[Lemma 3.6]{tao-noncommutative}
phrased using the language of control.

\begin{lemma}[Nonabelian Ruzsa covering lemma]
\label{ruz-cov} Suppose that $A$ and $B$ are finite sets in some ambient
group such that $B$ is symmetric and for which $|A \cdot B|, |B \cdot A|
\leq K|B|$ and $|B^2| \leq K|A|$. Then $A$ is $2K$-controlled by $B^2$.
\end{lemma}

\noindent\textit{Proof. } Take $X_1$ to be a maximal subset of $A$ with the
property that the translates $Bx$, $x \in X$, are all disjoint. It is easy
to see that $|X_1| \leq K$ and that $A \subseteq B^{-1} \cdot B \cdot X_1 =
B^2 \cdot X_1$. Similarly there is a set $X_2$ with $|X_2| \leq K$ such that 
$A \subseteq X_2 \cdot B^2$. Taking $X := X_1 \cup X_2$ we see that $A
\subseteq (B^2 \cdot X) \cap (X \cdot B^2)$, as required.\hfill{\usebox{\proofbox}}\vspace{11pt}

We supplement these last two results with the following lemma, which allows one to pass to finite index subgroups when describing the structure of approximate groups.
In view of Mal'cev's theorem it allows us to reduce the proof of Theorem \ref{mainthm2}
to the case of approximate subgroups of $\Upp_n(\mathbb{C})$.

\begin{lemma}
\label{finite-index-lem} Suppose that $G$ is a group and that $H\leq G$ is a
subgroup of $G$ with finite index $r$. Suppose that $A$ is a $K$-approximate subgroup of $G$. Then there is a $(rK)^{C}$-approximate
subgroup $S\subseteq H$ which $(rK)^{C}$-controls $A$.
\end{lemma}

\noindent \textit{Proof. } By the pigeonhole principle there is a coset $xH$
of $H$ such that $|A\cap xH|\geq |A|/r$. Writing $A^{\prime }:=A\cap xH$ and 
$B:=A^{\prime -1}A^{\prime }$, we have $B\subseteq H\cap A^{2}$. By Lemma \ref{ruz-cov} and Proposition \ref{approx-gp-basics} (ii) we see that $A$ is 
$(rK)^{C}$-controlled by $B^{2}$. By Proposition \ref{approx-gp-basics} (ii)
and (iv) there is a $(rK)^{C}$-approximate group $S\subseteq B^{\pm
6}\subseteq H$ which $(rK)^{C}$-controls $B^{2}$. By the transitivity of the
notion of control, $S$ also $(rK)^{C}$-controls $A$.\hfill {\usebox{\proofbox}}\vspace{11pt}

We saw in Proposition \ref{approx-gp-basics} that a $K$-approximate group $A$
satisfies the tripling condition $|A^{3}|\leq K^{2}|A|$. However there are
examples of sets with small tripling which are not approximate groups (an
example is a random subset of $\{1,\dots ,N\}$ of size $N/2$). For the
purposes of this paper it turns out to be convenient to work with the small
tripling condition rather than the approximate group condition. In the
remainder of this section we gather some preliminaries concerning sets with
small tripling (or even doubling).

\begin{lemma}[Iterated products]
\label{iterated-product} Suppose that $A$ is a subset of some ambient group
such that $|A^3| \leq K|A|$. Then we have $|A^{\pm n}| \leq K^{C_n}|A|$.
\end{lemma}

\noindent\textit{Proof. } This is not too hard to establish by repeated
applications of the Ruzsa triangle inequality. The details may be found in 
\cite[Lemma 3.4]{tao-noncommutative}.\hfill{\usebox{\proofbox}}\vspace{11pt}

\begin{lemma}[Almost uniform fibres]
\label{almost-uniform-fibres} Suppose that $G$ and $H$ are groups and that $\pi : G \rightarrow H$ is a homomorphism. Suppose that $A \subseteq G$ is a
set with $|A^2| \leq K|A|$. For each $x \in \pi(A)$ write $A_x$ for the
fibre $\pi^{-1}(x) \cap A$. Then $\max_{x \in \pi(A)} |A_x| \leq K \min_{x
\in \pi(A)}|A_x|$.
\end{lemma}

\noindent\textit{Proof. } Let $S \subseteq A$ be a set containing precisely
one element in each fibre $A_x$. Manifestly 
\begin{equation*}
|S| \min_{x \in \pi(A)} |A_x| \leq |A|.
\end{equation*}
On the other hand for any $x \in \pi(A)$ the set $A_x \cdot S$ has size $|A_x||S|$ and lies in $A \cdot A$, a set of size at most $K|A|$ by the
definition of what it means for $A$ to be a $K$-approximate group. Therefore 
\begin{equation*}
|S|\max_{x \in \pi(A)} |A_x| \leq K|A|.
\end{equation*}
Comparing these two inequalities yields the result.\hfill{\usebox{\proofbox}}
\vspace{11pt}

A fairly quick corollary of this is the next result, which examines the
extent to which small tripling is preserved under homomorphisms. Results
very similar to this one may be found in Section 7 of \cite{tao-noncommutative} and also in Lemma 7.3 of \cite{helfgott}.

\begin{corollary}[Small tripling preserved under homomorphisms]
\label{small-tripling-homs} Suppose that $G$ and $H$ are groups and that $\pi : G \rightarrow H$ is a homomorphism. Suppose that $A \subseteq G$ is a
set with $|A^3| \leq K|A|$. Then $|\pi(A)^3| \leq K^C|\pi(A)|$.
\end{corollary}

\noindent\textit{Proof. } Write $\approx$ to denote two quantities which are
equal up to some fixed power of $K$. We have $|A_x| \approx |A|/|\pi(A)|$
and $|(A^3)_x| \approx |A^3|/|\pi(A)^3|$ for all $x$ by the previous lemma
and Lemma \ref{iterated-product} (note that $\pi(A^3) = \pi(A)^3$, since $A$
is a homomorphism). Now it is clear that $A^3$ contains a translate of every
fibre of $A$, and so $\max_x |(A^3)_x| \geq \max_x |A_x|$. Combining these
facts with the assumption that $|A^3| \approx |A|$ gives $|\pi(A)^3| \approx
|\pi(A)|$, as required.\hfill{\usebox{\proofbox}}\vspace{11pt}

The behaviour of sets with small tripling under intersection with a subgroup 
$H \leq G$ is a little more subtle. Indeed suppose
that $G$ is large but finite and that $H$ is a small subgroup. If $S
\subseteq H$ is any set and if $A := S \cup (G \setminus H)$ then $A$ will
have tripling very close to 1 yet $S = A \cap H$ need not have any special
structure. This phenomenon is not particularly persistent, however, as the
following lemma shows. Lemmas almost identical to this one may be found in 
\cite{helfgott} and \cite{tao-solvable}.

\begin{lemma}[Sumsets and intersections]
\label{sumsets-intersections} Let $G$ be any group. Suppose that $A\subseteq
G$ is a finite symmetric set containing the identity with $|A^{3}|\leq K|A|$. Let $H\leq G$ be a subgroup. Then for each $n\geq 2$ we have 
\begin{equation*}
|A^{2}\cap H|\leq |A^{n}\cap H|\leq K^{C_{n}}|A^{2}\cap H|.
\end{equation*}
\end{lemma}

\noindent\textit{Proof. } Since $A$ contains the identity the lower bound is
obvious. To prove the upper bound list the right cosets of $H$ which have
nonempty intersection with $A$ as $Hx_1,\dots,Hx_k$, and set $A_{x_i} := A
\cap Hx_i$. Suppose without loss of generality that $|A_{x_1}|$ is the
largest or joint-largest of these sets. Of course it has cardinality at
least $|A|/k$. Then $A^2 \cap H$ contains $A_{x_1} A_{x_1}^{-1}$, and hence
also has cardinality at least $|A|/k$. Now the intersection of $A^{n+1}$
with each coset $H x_i$ has size at least $|A^n \cap H|$, and so $|A^n \cap
H| \leq |A^{n+1}|/k$. By Lemma \ref{iterated-product} this is at most $K^{C_n}|A|/k$. Comparing these inequalities gives the result. \hfill{\usebox{\proofbox}}\vspace{11pt}

Our final combinatorial preliminary is of a rather different type. It is a
beautiful ``sum-product'' result of Solymosi \cite{solymosi}.

\begin{lemma}[Sum-product over $\mathbb{C}$]
\label{sum-product} Suppose that $U,V,W$ are any three finite subsets of $%
\mathbb{C}$. Then we have 
\begin{equation*}
|U + V| |UW| \geq c|U|^{3/2}|V|^{1/2}|W|^{1/2}.
\end{equation*}
\end{lemma}

The use of sum-product phenomena in the study of approximate groups of
matrices is very natural, since matrix multiplication involves both addition
and multiplication of the entries. This observation has been exploited
before in the work of Helfgott \cite{helfgott-sl2,helfgott} and Chang \cite{chang1,chang2}.

\section{Proof of the main theorems}

By Lemma \ref{finite-index-lem} and Mal'cev's theorem it suffices to
establish Theorem \ref{mainthm2} for approximate subgroups of $\Upp_{n}(\mathbb{C})$. The argument proceeds by induction on $n$, and for the
purposes of establishing the inductive step it turns out to be advantageous
to work with the following (slightly more precise) variant of Theorem \ref{mainthm2}.

\begin{mainthm2-again}
Let $K \geq 1$. Suppose that $G = \Upp_n(\mathbb{C})$ and that $A
\subseteq G$ is a set with $|A^3| \leq K|A|$. Then there is some set $A^{\prime }\subseteq A$ with $|A^{\prime}| \geq K^{-C}|A|$ which is
contained in a left coset of a nilpotent subgroup of $\Upp_n(\mathbb{C})$ of step at most $n-1$.
\end{mainthm2-again}

\emph{Deduction of Theorem \ref{mainthm2} for approximate subgroups of $\Upp_{n}(\mathbb{C})$.} Let $A\subseteq \Upp_{n}(\mathbb{C})$ be
a $K$-approximate group. Then by Proposition \ref{approx-gp-basics} (i) we
have $|A^{3}|\leq K^{2}|A|$, and so we may apply Theorem \ref{mainthm2}'
(with $K$ replaced by $K^{2}$). It follows that there is some set $A^{\prime
}\subseteq A$ with $|A^{\prime }|\geq K^{-C}|A|$ which is contained in $xH$,
some left coset of a nilpotent group $H\leq \Upp_{n}(\mathbb{C})$ of
step at most $n-1$. The set $S=A^{\prime -1}A^{\prime }$ then generates a
nilpotent subgroup of $\Upp_{n}(\mathbb{C})$, and by Lemmas \ref{ruz-cov} and \ref{iterated-product} the set $S^{2}$ is symmetric and $K^{C}$
-controls $A$. By the second sentence of Proposition \ref{approx-gp-basics}
(v) and Lemma \ref{iterated-product} the set $S^{6}$ is a $K^{C}$-approximate group which $K^{C}$-controls $A$. This concludes the proof of
Theorem \ref{mainthm2}.\hfill {\usebox{\proofbox}}\vspace{11pt}

It remains, then, to establish Theorem \ref{mainthm2}'. As we said this
proceeds by induction, the case $n = 1$ being trivial since $\Upp_1(\mathbb{C}) \cong \mathbb{C}$ is abelian. The key to the inductive step is
the observation that there are two different homomorphisms $\pi,\pi^{\prime
}: \Upp_n(\mathbb{C}) \rightarrow \Upp_{n-1}(\mathbb{C})$
defined by 
\begin{equation*}
\pi\big((x_{ij})_{1 \leq i \leq j \leq n}\big) = (x_{ij})_{1 \leq i \leq j
\leq n-1} \qquad \mbox{and} \qquad \pi^{\prime }\big((x_{ij})_{1 \leq i \leq
j \leq n}\big) = (x_{i+1,j+1})_{1 \leq i \leq j \leq n-1}.
\end{equation*}
Let us take a set $A \subseteq \Upp_n(\mathbb{C})$ with $|A^3| \leq
K|A|$ and apply the two homomorphisms $\pi$ and $\pi^{\prime }$ to it. By
Lemma \ref{small-tripling-homs} we have $|\pi(A)^3| \leq K^4|\pi(A)|$ and
hence by inductive hypothesis there is some set $A^{\prime }\subseteq A$
with $|\pi(A^{\prime})| \geq K^{-C}|\pi(A)|$ such that $\pi(A^{\prime})$ is
contained in a left coset of a nilpotent subgroup of $\Upp_{n-1}(\mathbb{C})$. It follows from Lemma \ref{almost-uniform-fibres} that $%
|A^{\prime}| \geq K^{-C}|A|$. We now apply $\pi^{\prime }$ to \emph{this}
set, obtaining a further subset $A^{\prime \prime }\subseteq A^{\prime }$
with $|A^{\prime \prime}| \geq K^{-C}|A|$ and such that $\pi^{\prime
}(A^{\prime \prime })$ is contained in a left coset of a nilpotent subgroup
of $\Upp_{n-1}(\mathbb{C})$. For ease of notation we drop the double
dash and assume that both $\pi(A^{\prime})$ and $\pi^{\prime }(A^{\prime })$
lie in left cosets of nilpotent subgroups of $\Upp_{n-1}(\mathbb{C})$
of step at most $n-1$. The sets $\pi(A^{\prime} A^{\prime -1})$ and $\pi^{\prime }(A^{\prime} A^{\prime -1})$ are then contained in nilpotent
subgroups of $\Upp_{n-1}(\mathbb{C})$ of step at most $n-1$. Write $B
:= A^{\prime} A^{\prime - 1}$; thus $B$ is symmetric and by Lemma \ref{iterated-product} we have 
\begin{equation}  \label{B-tripling}
|B^3| \leq |A^6| \leq K^{C}|A| \leq K^{C_d}|B|.
\end{equation}
Our conclusion implies that for any $b_1,\dots b_{n} \in B$ the commutator 
\begin{equation*}
[b_1,[b_2,[b_3,\dots [b_{n-1},b_n]\dots]]]
\end{equation*}
lies in $\ker \pi \cap \ker \pi^{\prime }$, which is none other than the
normal subgroup $H \lhd \Upp_n(\mathbb{C})$ consisting of matrices of
the form $m_{\lambda} = \left(\begin{smallmatrix}
1 & 0 & \cdots & \lambda \\ 
0 & 1 & \cdots & 0 \\ 
\vdots & \vdots & \ddots & \vdots \\ 
0 & 0 & \cdots & 1\end{smallmatrix}\right)$, $\lambda \in \mathbb{C}$. Note that $m_{\lambda}m_{\mu} =
m_{\lambda + \mu}$.

We will make use of the following well-known lemma:

\begin{lemma}
\label{nil-2} Suppose that $G$ is a group generated by finitely many
elements $x_1,\dots,x_k$, and that every $(s+1)$-fold nested commutator $[x_{i_1},[x_{i_2},[x_{i_3}, \dots [x_{i_s}, x_{i_{s+1}}]\dots]]]$ is equal
to the identity. Then $G$ is $s$-step nilpotent.
\end{lemma}

\noindent\textit{Proof. } This is \cite[Lemma 8.17]{Rag}. \hfill{\usebox{\proofbox}}\vspace{11pt}

Let $N := 3\cdot 2^{n-1} - 2$, the length of the commutator $[b_1,[b_2,[b_3,\dots [b_{n-1},b_n]\dots]]]$ as a word in the $b_i$. Now if $B^N \cap H$ contains only the identity matrix then every $n$-fold nested
commutator $[b_1,[b_2,[b_3,\dots [b_{n-1},b_{n}]\dots]]]$, where $b_1,\dots,b_{n} \in B$, is equal to the identity. By Lemma \ref{nil-2} it
follows that $B$ generates an $(n-1)$-step nilpotent group.

Otherwise, $B^N \cap H$ contains at least one matrix $m_{\lambda}$ with $\lambda \neq 0$. Note that $B^N$ also contains the identity matrix. We now
divide into two cases.

Let $D = K^{\gamma_n}$, where $\gamma_n > 0$ is to be specified later.
Suppose first that there are at least $D$ different ratios $T =
\{t_1,\dots,t_D\}$ amongst the values of $x_{11}/x_{nn}$ occurring in $B^{N}$
. For each $k = 1,\dots,D$ choose some $y_k = (x^{(k)}_{ij}) _{1 \leq i\leq j \leq n}\in B^{N}$ with $x^{(k)}_{11}/x^{(k)}_{nn} = t_k$. An easy computation gives $y_k m_{\lambda} y_k^{-1} = m_{t_k\lambda}$. Thus if the set $S \subseteq \mathbb{R}$ is defined by 
\begin{equation*}
B^{N} \cap H = \{m_{\lambda} : \lambda \in S\}
\end{equation*}
then $B^{3N} \cap H$ contains the elements $\{m_{\lambda} : \lambda \in TS\}$
and also the elements $\{m_{\lambda} : \lambda \in S + S\}$. However by
Solymosi's sum-product bound (Lemma \ref{sum-product}) we have 
\begin{equation*}
|S + S| \cdot |ST| \geq c |S| |T|^{1/2}.
\end{equation*}
Thus 
\begin{equation*}
|B^{3N} \cap H| \geq c \sqrt{D} |B^{N} \cap H|.
\end{equation*}
If $\gamma_n$ is sufficiently large then, together with Lemma \ref{B-tripling}, this contradicts Lemma \ref{sumsets-intersections}.

Such a set $T$ cannot, therefore, exist. \emph{A fortiori} there are no more
than $D$ different ratios $x_{11}/x_{nn}$ occurring amongst the elements $(x_{ij})$ of $A^{\prime }$, and so there is a set $C \subseteq A^{\prime}
A^{\prime -1} = B$, of size at least $|B|/D \gg K^{-C_n}|A|$, with the
property that $x_{11} /x_{nn} = 1$ for all $(x_{ij}) \in C$. Then every
element of $C$ commutes with $H$ and so 
\begin{equation*}
[c,[b_1,[b_2,\dots [b_{n-1},b_{n}]\dots ]]] = \id
\end{equation*}
whenever $b_1,\dots,b_n \in B$ and $c \in C$, and hence in particular when
all of $b_1,\dots,b_n,c$ lie in $C$. It follows from Lemma \ref{nil-2} that $%
C$ is contained in a nilpotent subgroup of $\GL_n(\mathbb{C})$.

We have examined two cases, and in both of them have located a set $X
\subseteq A^{\prime} A^{\prime -1}$ with $|X| \geq K^{-C_n}|A|$ which is
contained in a nilpotent subgroup of $\Upp_n(\mathbb{C})$ of step at
most $n$. Remember what it is that we are trying to prove: that a
substantial fraction of $A$ lies inside some \emph{coset} of a nilpotent
subgroup of $\GL_n(\mathbb{C})$. To establish this it suffices to show
that there is $x \in \GL_n(\mathbb{C})$ such that $|A \cap xC| \gg
K^{-O_n(1)}|A|$. Write $r(x) = |A \cap xC|$, which is also the number of
representations of $x$ as $ac^{-1}$ with $a \in A, c \in C$. Then 
\begin{equation*}
\sum_x r(x) = |A||C| \gg K^{-O_n(1)}|A|^2.
\end{equation*}
However $r(x) = 0$ unless $x \in AC^{-1}$, a set which is contained in $A^3$
and hence has cardinality at most $K|A|$. It follows that there is at least
one $x$ for which $r(x) \gg K^{-O_n(1)}|A|$, as required. This completes the
proof of the inductive step in Theorem \ref{mainthm2}', and hence the proof
of Theorem \ref{mainthm2}.\hfill{\usebox{\proofbox}}\vspace{11pt}

It remains to deduce Corollary \ref{main-cor}. To do this we need the
following proposition concerning the structure of nilpotent subgroups of $\Upp_n(\C)$. 
Here $\mathbb{T}^{n}:=\mathbb{R}^{n}/\mathbb{Z}^{n}$ is the $n$-dimensional
torus.

\begin{proposition}
\label{prop-a3} Suppose that $G$ is a nilpotent subgroup of $\Upp_{n}(\mathbb{C})$. Then $G$ embeds into $\mathbb{T}^{n}\times \Gamma $, where $\Gamma $ is a simply-connected nilpotent Lie group of dimension at most that
of $\Upp_{n}(\mathbb{C})$ \textup{(}and in particular is torsion-free\textup{)}.
\end{proposition}

Although this proposition is certainly well-known to experts, we do not know of a particularly convenient source for it. The account which follows assumes a certain amount of background in the theory of complex algebraic groups such as may be found in the book of Humphreys \cite{humphreys}. A more self-contained discussion based on the treatment of Wehrfritz \cite{wehrfritz} may be found on the second author's website \cite{green-notes}.\vspace{11pt}

\noindent\textit{Proof. } As $G$ embeds in its Zariski-closure, and the Zariski-closure of a nilpotent group is still nilpotent, there is no loss of generality in assuming that $G$ is Zariski-closed. We shall show that $G$ is a
direct product of two Zariski-closed subgroups $G_{s}$ and $G_{u}$ where $G_{s}$ is isomorphic to a subgroup of $(\mathbb{C}^{\times })^{n}$ and $G_{u}
$ is connected, simply connected and nilpotent. This is enough to prove the
proposition, since $(\mathbb{C}^{\times })^{n}$ is $\mathbb{T}^{n}\times 
\mathbb{R}^{n}$ as a real Lie group. When $G$ is Zariski-connected, this
assertion is Proposition 19.2 in Humphreys' book \cite{humphreys}. We recall the argument and explain why it continues to hold in the
non-connected case. 

Let $G_{u}$ be the set of unipotent elements in $G$. Clearly $G_{u}=G\cap\Uni(\mathbb{C})$, where $\Uni(\mathbb{C})$ is the subgroup of $\Upp(\mathbb{C})
$ consisting of matrices with all eigenvalues equal to $1$. First we remark that $G_{u}$ is connected (in either the topological or Zariski- sense, these being equivalent for complex algebraic groups). To see this note that if $\{g^{n}\}_{n\in \mathbb{Z}}$ lies in $G_{u}$ so does $\{g^{t}\}_{t\in \mathbb{C}}$, since $g^n$ is a polynomial function of $n$, and so $G_u$ is path-connected. Now it is known (see for example \cite{bourbaki}) that the exponential map on a connected, simply-connected Lie group $H$ such as $\Uni(\mathbb{C})$ is a diffeomorphism and that it induces a correspondence between closed connected subgroups of $H$ and Lie subalgebras of $\mathfrak{h} := \log H$. These last are certainly vector subspaces of $\mathfrak{h}$ and hence are simply-connected, and hence so are all closed connected subgroups of $H$. In particular, $G_u$ is simply-connected.

Let $G_{s}$ be the set of semisimple elements in $G$. We claim that $G_{s}$ is a subgroup and that it commutes pointwise with $G_{u}$. Once this is established it follows that $G=G_{s}\times G_{u}$, and furthermore the natural map $\Upp(\mathbb{C)}\rightarrow (\mathbb{C}^{\times })^{n}$ restricts to an isomorphism on $G_{s}$, thereby concluding the proof of the proposition. In fact it is enough to show that $G_{s}$ commutes pointwise with $G_{u}$. To see this, let $x,y \in G_s$ be arbitrary and consider the commutator $z := xyx^{-1}y^{-1}$, which lies in $G_u$. Then we have
\[ xyx^{-1} = z \cdot y = y \cdot z\] and
\[ xyx^{-1} = \id \cdot (xyx^{-1}) = (xyx^{-1})\cdot \id ,\] two decompositions of $xyx^{-1}$ into commuting semisimple and unipotent parts. By the uniqueness of Jordan decomposition it follows that $xyx^{-1} = y$, and so $xy = yx$. But the product of two commuting semisimple transformations is semisimple, whence $xy \in G_s$.

To confirm the claim, then, it remains to establish that $G_s$ commutes with $G_u$. To this end we quote \cite[Theorem 18.3 (c)]{humphreys}, which tells us that for each fixed $y\in
G_{s}$ the map $\phi_y(x) := xyx^{-1}y^{-1}$ is a bijection from $\phi_y(G_u)$ to itself.  Thus we can write each element of $\phi_y(G_u)$ as a commutator of arbitrarily high order involving elements of $G$. Since $G$ is nilpotent, we are forced to conclude that $\phi_y(G_u) = \{\id\}$ and hence that $G_s$ does indeed commute with $G_u$.\endproof

\emph{Deduction of Corollary \ref{main-cor}.} The argument is very close in
spirit to the \textquotedblleft rectification principle\textquotedblright\
developed in \cite{bilu-lev-ruzsa}. Let $A$ be a $K$-approximate subgroup of
some solvable subgroup of $\GL_{n}(\mathbb{C})$. By Theorem \ref{mainthm2} there is a $K^{C_{n}}$-approximate subgroup $B\subseteq \langle
A\rangle $ generating a nilpotent subgroup of $\Upp_{n}(\mathbb{C})$
which $K^{C_{n}}$-controls $A$. By Proposition \ref{prop-a3} this nilpotent
subgroup embeds into $\mathbb{T}^{n}\times \Gamma $, where $\Gamma $ is a
simply-connected nilpotent Lie group. It therefore suffices to show that if $S\subseteq \mathbb{T}^{n}\times \Gamma $ is a $K$-approximate group then it
is $e^{K^{C_{n}}}$-controlled by a nilpotent progression in $\langle
S\rangle $. By a trivial pigeon-hole argument there is a set $S^{\prime
}\subseteq S$, $|S^{\prime }|\geq 6^{-n}|S|$, which is entirely contained in
some box $\mathbf{x}+[-\frac{1}{12},\frac{1}{12}]^{n}\subseteq \mathbb{T}^{n}
$. Let $S^{\prime \prime }=S^{\prime }S^{\prime -1}.$ We have $S^{\prime
\prime }\subseteq S^{2}$, $|S^{\prime \prime }|\geq 6^{-n}|S|$ and $S^{\prime \prime }\subseteq \lbrack -\frac{1}{6},\frac{1}{6}]^{n}\subseteq 
\mathbb{T}^{n}$. Consider the obvious lift $\psi :S^{\prime ^{\prime
}}\rightarrow \lbrack -\frac{1}{6},\frac{1}{6}]^{n}\times \Gamma \subseteq 
\mathbb{R}^{n}\times \Gamma $ to the universal cover. It is clear that $|\psi (S^{\prime \prime })^{\pm 3}|=|S^{\prime \prime \pm 3}|$, and so by
Proposition \ref{approx-gp-basics} $\psi (S^{\prime \prime })^{\pm 3}$ is a $K^{C_{n}}$-approximate group which $K^{C_{n}}$-controls $\psi (S^{\prime
\prime })$. By the main result of \cite{breuillard-green-nil} it follows that there is a nilpotent progression $
P\subseteq \langle \psi (S^{\prime \prime })\rangle $ of dimension at most $
K^{C_{n}}$ which $e^{K^{C_{n}}}$-controls $\psi (S^{\prime \prime })^{\pm 3}$
. The image of this progression under the canonical projection back down to $
\mathbb{T}^{n}\times \Gamma $ is a nilpotent progression $Q$ of the same
dimension which is contained in $\langle S\rangle $, and which $e^{K^{C_{n}}}
$-controls the set $S^{\prime \prime \pm 3}$. But by Proposition \ref{approx-gp-basics} (v) $S^{\prime \prime \pm 3}$ is a $K^{C_{n}}$
-approximate group which $K^{C_{n}}$-controls $S.$ Hence $S$ is $e^{K^{C_{n}}}$-controlled by $Q.$ \endproof

\end{document}